\documentclass[a4paper,12pt]{amsart}
\usepackage{amsmath,amssymb}
\def\C{\Bbb C}
\def\CC{\mathcal C}
\def\D{\Bbb D}
\def\DD{\Delta}
\def\O{\mathcal O}
\def\d{\delta}
\def\e{\eta}
\def\z{\zeta}
\def\eps{\varepsilon}
\def\g{\gamma}
\def\k{\kappa}

\def\vp{\varphi}
\def\ov{\overline}
\def\bd{\partial}
\def\F{\mathcal F}
\def\DD{\mathcal D}
\def\ds{\displaystyle}
\newtheorem{thm}{Theorem}
\newtheorem{prop}[thm]{Proposition}
\newtheorem{cor}[thm]{Corollary}

\title[Comparison and localization of invariant functions]
{Comparison and localization of invariant functions
on strongly pseudoconvex domains}

\author{Nikolai Nikolov}
\address{N. Nikolov\\Institute of Mathematics and Informatics\\Bulgarian Academy
of Sciences\\Acad. G. Bonchev 8, 1113 Sofia, Bulgaria
\vspace{1mm}
\newline Faculty of Information Sciences\\
State University of Library Studies and Information Technologies\\
Shipchenski prohod 69A, 1574 Sofia,
Bulgaria}

\email{nik@math.bas.bg}

\thanks{The author was partially supported by the Bulgarian National Science Fund,
Ministry of Education and Science of Bulgaria under contract KP-06-Rila/2.}


\begin{document}


\begin{abstract} Comparison and localization results for the Lempert function,
the Carath\'eodory distance and their infinitesimal forms on strongly
pseudoconvex domains are obtained. Related results for visible and strongly complete
domains are proved.
\end{abstract}

\maketitle

\section{Introduction and main results}

Let $D$ be a domain in $\C^n,$ $z,w\in D$ and $X\in\C^n.$
The Kobayashi distance $k_D$ is the largest pseudodistance
not exceeding the Lempert function
$$l_D(z,w)=\inf\{\tanh^{-1}|\alpha|:\exists\vp\in\O(\D,D)
\hbox{ with }\varphi(0)=z,\varphi(\alpha)=w\},$$ where $\D$ is the
unit disc and $\tanh^{-1}t=\frac12\log\frac{1+t}{1-t}.$
Note that $k_D$ is the integrated form of the Kobayashi metric
$$\k_D(z;X)=\inf\{|\alpha|:\exists\vp\in\O(\D,D)\hbox{ with }
\varphi(0)=z,\alpha\varphi'(0)=X\}.$$

The Carath\'eodory distance and metric are given by
$$c_D(z,w)=\sup\{\tanh^{-1}|f(w)|:f\in\O(D,\D)
\hbox{ with }f(z)=0\},$$
$$\g_D(z;X)=\sup\{|f'(z)X|:f\in\O(D,\D)\},$$
respectively. Note that $c_D\le k_D\le l_D$ and $\g_D\le\k_D.$

Lempert's theorem implies that (see \cite[Theorem 1]{L2}) $c_D=l_D$ and, therefore,
$\g_D=\k_D$ on any convex domain.

Let now $D$ be strongly pseudoconvex. Then
\begin{equation}\label{root}l_D(z,w)-c_D(z,w)\le C|z-w|^{1/2},
\end{equation}
\begin{equation}\label{rt}\k_D(z;X)-\g_D(z;X)\le C|X|.
\end{equation}
by \cite[Theorem 1.6]{NT2} and \cite[Theorem 3.2]{Fu}, resp.

The first aim of this paper is to obtain better estimates than \eqref{root} and \eqref{rt}.
Set $\d_D(z)=\mbox{dist}(z,\bd D),$ $\ds h(x)=\frac{x(1+x)}{\log(1+x)},\ x>0,$ and
$$f_D(z,w)=\d_D(z)\d_D(w)h\left(\frac{|z-w|}{\d_D(z)^{1/2}\d_D(w)^{1/2}}\right),$$
$$g_D(z,w)=|z-w|(|z-w|+\d_D(z)^{1/2}\d_D(w)^{1/2}).$$

\begin{thm}\label{global} Let $D$ be a strongly pseudoconvex domain. Then there
exists $C>0$ such that for $z,w\in D,$ $z\neq w,$ and $X\in\C^n$ one has that
\begin{equation}\label{zero}\frac{l_D(z,w)}{c_D(z,w)}\le 1+Cf_D(z,w),
\end{equation}
\begin{equation}\label{one}
l_D(z,w)-c_D(z,w)\le Cg_D(z,w),
\end{equation}
\begin{equation}\label{two}
\k_D(z;X)\le(1+C\d^2_D(z))\g_D(z;X)\le\g_D(z;X)+C\d_D(z)|X|.
\end{equation}
\end{thm}

The main point in the proof of Theorem \ref{global} is the following localization result.

\begin{prop}\label{local} Let $D$ be a bounded pseudoconvex domain such that either $\ov D$
is a Stein compact or $\bd D$ is $\CC^\infty$-smooth. Let $p\in\bd D$ be a strongly
pseudoconvex boundary point. Then there exists a neighborhood $U_0$ of $p$ such that for any
neighborhoods $V\Subset U\subset U_0$ of $p$ one may find $C>0$ such that
for $z,w\in D\cap V,$ $z\neq w,$ and $X\in\C^n$ one has that
\begin{equation}\label{quo}
\frac{k_{D\cap U}(z,w)}{c_D(z,w)}\le 1+Cf_D(z,w),
\end{equation}
\begin{equation}\label{lip}
k_{D\cap U}(z,w)-c_D(z,w)\le Cg_D(z,w),
\end{equation}
\begin{equation}\label{met}
\k_{D\cap U}(z;X)\le(1+C\d^2_D(z))\g_D(z;X)\le\g_D(z;X)+C\d_D(z)|X|.
\end{equation}
\end{prop}

Here and what follows $D\cap U$ is assumed to be connected.

These estimates are optimal on $\D$ and, hence, on the unit ball in
$\C^n.$ More precisely, if $D=\D,$ $r\in(0,2),$ $U_0=\D(1,r),$ and $z,w\in\Bbb R,$
then an opposite inequality to \eqref{lip} holds (see Remark (a) after the proof of Proposition
\ref{convex}). This easily implies the same for \eqref{quo} and \eqref{met}. So there are no better
upper bounds involving only $|z-w|,$ $\d_D(z)$ and $\d_D(w).$

The inequality \eqref{lip} is also valid (see \cite[Proposition 4.2]{NT2}) for any planar domain $D$ with
Dini-smooth boundary near $p\in\bd D.$ (For $D\subset\Bbb R^d,$ this means that the inner unit normal
vector to $\CC^1$-smooth $\bd D$ near $p$ is a Dini-continuous function.)
One may show the same for \eqref{quo} and \eqref{met}.

The inequalities \eqref{zero}, \eqref{quo}, and \eqref{met} can be considered as quantitative versions
of \cite[Theorem 1, Proposition 3, and Proposition 1]{Ven} which say that
the respective quotients tend to 1 when $z\to\bd D/z,w\to p.$

The inequality \eqref{met} is covered by \cite[Theorem 1.2]{Wol} when $z$ lies on the inward normal
to $\CC^4$-smooth $\bd D$ at $p$ and $X$ is parallel to this normal.

Recall that $\g_D\le\k_D\le\k_{D\cap U}$ and $c_D\le k_D\le k_{D\cap U},$ in general.
On the other hand, it follows by \cite[Theorem 1]{H1} and \cite[Proposition 2.5]{BFW}, resp., that
if $D$ is strongly pseudoconvex and $\bd D$ is $\CC^3$-smooth near $p,$ then one may
find $\eps>0$ such that $\k_{D\cap U}(z;X)=\g_D(z;X)$
and $k_{D\cap U}(z,w)=c_D(z,w)$ if $|X_N|<\eps|X|,$ $|z-p|<\eps,$ $|w-p|<\eps,$ and $|(z-w)_N|<\eps|z-w|,$
where $|Y_N|$ is the projection of $Y$ on the complex normal to $\bd D$ at the closest point to $z$
(a related result can be found in \cite[Theorem 2]{H2}). A recent observation by X. Huang \cite[p.~1838]{H3}
shows that $\CC^{2,\alpha}$-smoothness is enough.

So, if $D$ is strongly pseudoconvex and $\bd D$ is $\CC^{2,\alpha}$-smooth, then there exists $\eps>0$
such that $l_D(z,w)=c_D(z,w)$ and $\k_{D\cap U}(z;X)=\g_D(z;X)$ if $\d_D(z)<\eps,$ $\d_D(w)<\eps,$
$|(z-w)_N|<\eps|z-w|,$ and $|X_N|<\eps|X|.$

Note also that, assuming only strongly pseudoconvexity near $p,$ weaker inequalities than \eqref{lip}
and \eqref{met} are known (see \cite[Theorem 1.2 (1.11)]{NT2} and \cite[Theorem 2.1]{FR}, resp.):
$$k_{D\cap U}(z,w)-k_D(z,w)\le C|z-w|^{1/2},\quad z,w\in D\cap V,$$
$$\k_{D\cap U}(z;X)\le(1+C\d_D(z))\k_D(z;X),\quad z\in D\cap V,\ X\in\C^n.$$

The proof of Proposition \ref{local} uses the basic idea of X. Huang \cite{H1} to reduce the
argument to a localized version of the Forn{\ae}ss embeding  theorem (see \cite[Proposition 1]{For}
and \cite[Lemma 5]{H1}) and Proposition \ref{convex} below.

\section{Visibility and strong completeness}

\subsection{Visible domains} Let $D$ be a bounded complete domain (w.r.t. $k_D$).
Then any two points can be joined by a real geodesic,
i.e. a curve $\psi:[a,b]\to D$ such that $k_D(\psi(s),\psi(t))=|s-t|$ for any $s,t\in[a,b].$
We say that $D$ is (weakly) \emph{visible} if for any $\eps>0$ one can find $K\Subset D$ that intersects
any real geodesic $\psi$ with $\mbox{diam}(\psi)>\eps.$ One may easily see that this is equivalent to
the usual definition of visibility: for any $p,q\in\partial D,$ $p\neq q,$ there exist neighborhoods
$U$ of $p$ and $V$ of $q,$ and $K\Subset D$ such that any real geodesic $\psi$ that meets $U$ and $V$
must meet $K,$ too.

Denote by
$$(z|w)_o=k_D(z,o)+k_D(w,o)-k_D(z,w).$$
the doubled Gromov product. It turns out (see e.g. \cite[Proposition 2.5]{BNT}) that visibility is equivalent to
\begin{equation}\label{grom}
\limsup_{z\to p,w\to q}(z|w)_o<\infty
\end{equation}
for any $p,q\in\bd D,$ $p\neq q$ and some (hence any) $o\in D.$

Note that any bounded convex domain $D$ is visible if
\begin{equation}\label{gold}
\hat\d_D(z)\le f(\d_D(z))\mbox{ near }\bd D,
\end{equation}
where $\int_0^\eps\frac{f(x)}{x}<\infty$ and $\hat\d_D(z)$ denotes the maximal radius
of complex affine discs in $D$ centered at $z$ (see \cite[Theorem 1.4]{BZ1} for a more general result).
The above condition requires some control on complex flatness of $\bd D.$ For example,
it is satisfied if $f(x)=Cx^{1/m},$ $m>0;$ then $D$ is called $m$-convex (see \cite[Definition 2.7]{Mer}).
This notion extends the notion of finite type to the non-smooth case. So any strongly convex domain
is 2-convex and hence visible.

One may have visibility even if $\bd D$ is arbitrarily complex-flat. Indeed, if $D$ is a bounded domain
with Dini-smooth boundary, then (see \cite[Corollary 8]{NA})
\begin{equation}\label{dini}
k_D(z,w)\le\log\left(1+\frac{C|z-w|}{\d_D(z)^{1/2}\d_D(w)^{1/2}}\right).
\end{equation}
On the other hand, if $D$ is bounded and convex, then (see e.g.~\cite[(2.4)]{BNT})
\begin{equation}\label{low}
k_D(z,w)\ge\frac{1}{2}\log\left|\frac{\d_D(z)}{\d_D(w)}\right|.
\end{equation}
Combining \eqref{grom}, \eqref{dini} and
\eqref{low} implies that visibility of a bounded convex domain $D$ with Dini-smooth boundary means that
for any $p,q\in\bd D,$ $p\neq q,$ and any $z\in D$ near $p,$ $w\in D$ near $q,$ one has that
\begin{equation}\label{vis}
k_D(z,w)\ge\frac{1}{2}\log\frac{1}{\d_D(z)}+\frac{1}{2}\log\frac{1}{\d_D(w)}-C.
\end{equation}
The last condition is satisfied if, for example, $D$ is strictly $\C$-convex,
i.e. $\bd D\cap T^\C_p(\bd D)=\{p\}$ for any $p\in\bd D$
(see \cite[Subsection 2.4.1 and Lemma 4.5]{Zim}).

To prove Proposition \ref{local}, we will use its general version in the convex setting.

\begin{prop}\label{convex} Let $D$ be a bounded convex visible domain with Dini-smooth boundary.
Then for any $p\in\bd D$ and any neighborhoods $V\Subset U$ of $p$ there exists $C>0$
such that \eqref{quo}, \eqref{lip}, and \eqref{met} hold for $z,w\in D\cap V,$ $z\neq w,$ and $X\in\C^n$.
\end{prop}

For other localization results see e.g. \cite[Theorem 1.4 and Proposition 6.15]{BNT}.

\subsection{Strongly complete domains}
The proof of Proposition \ref{convex} involves the boundary behavior of complex geodesics,
i.e. maps $\vp\in\O(\D,D)$ such that $k_\D(\z,\e)=k_D(\vp(\z),\vp(\e))$ for any $\z,\e\in\D.$

Some results in this direction can be found, for example, in \cite{L1,L3,Mer,H1,JP}.

Before we state our result, we need the following definition: a bounded complete domain $D$
is said to be \emph{strongly complete} if for any sequence $(\psi_j)$ of real geodesics
joining $z_j,w_j\to p\in\partial D$ one has that
$$k_D(o,\psi_j)\to\infty$$
for some (hence any) $o\in D.$

In \cite{BZ2}, such a domain is said to to have well-behaved geodesics.

Note that the following counterpart to \eqref{grom}: if
$$\lim_{z,w\to p}(z|w)_o=\infty$$
for any $p\in\bd D,$ then $D$ is strongly complete.

This follows by the inequality $2k_D(o,u)\ge (z|w)_o$ for any $u$ on a
real real geodesic $z$ and $w.$

Observe also that \eqref{dini} and \eqref{low} implies that any bounded convex domain with
Dini-smooth boundary is strongly complete (see also \cite[Lemma 4.1]{LW}).

The next proposition is inspired by \cite[Theorem 1.6]{BZ2} (see also \cite[Theorem 2.18]{Zim},
\cite[Proposition 2.3]{BFW}, and \cite[Corollary 4.2 and Remark 4.3]{LW}).

\begin{prop}\label{equi} Let $G$ be a strongly complete domain and $D$ a visible domain. Let
$\tau\ge 0,$ $K\Subset G,$ and $\ds\lim_{x\to 0}g(x)=0.$ Then the family $\F$ (possibly empty) of
Kobayashi isometries $f:G\to D$ such that \footnote{In fact, $\sup_G=\max_G$ since $G$ is complete.}
$\sup_G(\d_D\circ f)\le g(\sup_K(\d_D\circ f))$ is uniformly equicontinuous
(with respect to the Euclidean distance).
In particular, any Kobayashi isometry $f:G\to D$ extends continuously on $\ov G$
and $f(\bd G)\subset\bd D.$
\end{prop}

Taking $G=\D,$ Proposition \ref{equi} covers the case of complex geodesics.

Note also the main result in \cite{Mai}, Theorem 1.3, is exactly the second statement above
in the particular case when $G$ and $D$ are bounded convex domain with Dini-smooth boundaries,
and $D$ is strictly $\C$-convex.

Recall now that if $D$ is a bounded convex domain, then for any $z,w\in D,$ $z\neq w,$ and any
$X\in\C^n,$ $X\neq 0,$ there exist extremal maps for $l_D(z,w)$ and $\k(z;X)$ and any such map
turns out to be a complex geodesic (see \cite{RW}).

We call $\bd D$ disc-free/affine disc-free if any holomorphic/complex affine disc in
$\bd D$ is constant. When $D$ is convex, both notions coincide (see e.g. \cite[Theorem 1.1]{FS}).

Proposition \ref{equi} suggests the next characterization of visibility.

\begin{prop}\label{geod} A bounded convex domain $D$ is visible if and only if $\bd D$
is disc-free and the family $\F$ of complex geodesics $\vp:\D\to D$ such that
$(\ast)$ $\d_D \circ \vp(0)=\max(\d_D\circ\vp)$ is uniformly equicontinuous.\footnote{Any geodesic
can be reparameterized to satisfy $(\ast)$.}
\end{prop}

The second part of the proof of this proposition easily leads to the following consequence.

\begin{cor}\label{cor} Let $D$ be a visible convex domain. Then for any $p,q\in\ov D,$
$p\neq q,$ there exists a complex geodesic $\vp$ that extends continuously on $\ov\D$ and
$p,q\in\vp(\ov\D).$
\end{cor}

\section{Proofs}

\noindent{\it Proof of Proposition \ref{equi}.} Assume that the family $\F$ is not uniformly
equicontinuous. Then we may find sequences $z_j,w_j\to b\in\bd G,$ $(o_j)\subset\ov K,$
and $(f_j)\subset\F$ such that $z_j'=f_j(z_j)\to p,$ $w_j'=f_j(w_j)\to q,$ $o_j'=f_j(o_j)\to s,$
where $p\neq q\neq s\in\ov D$ and $\d_D(o_j')=\sup_K(\d_D\circ f_j).$

Assume that $s\in\bd D.$ Then $f_j\to\bd D$ uniformly on $G.$ Let $\psi_j$ be a real geodesic
joining $o_j$ and $w_j.$ Then $f_j\circ\psi_j\to\bd D$ is a real geodesic joining $o_j'\to s$ and
$w_j'\to q.$ Since $s\neq q,$ we get a contradiction with the visibility of $D.$

Let now $s\in D.$ Let $\vp_j$ be a real geodesic joining $z_j$ and $w_j.$ Then $f_j\circ\vp_j\to\bd D$
is a real geodesic joining $z_j'\to p$ and $w_j'\to q,$ and
$$k_D(o_j',f_j\circ\vp_j)=k_G(o_j,\vp_j)\to\infty$$
because $s\in D$ and $G$ is strongly complete. Hence $f_j\circ\vp_j\to\partial D.$
Since $p\neq q,$ this contradicts to the visibility of $D.$

Finally, the the uniform continuity of $f$ implies that $f$ extends continuously on $\ov G,$
and $f(\bd G)\subset\bd D$ follows by completeness of $G.$\qed
\smallskip

\noindent{\it Remark.} This proof implies that any Kobayashi isometry $f:G\to D$ extends continuously
on $\ov G$ if the visibility of $D$ is replaced by the weaker assumption that $D$ is bounded and
satisfies \eqref{grom}.
\medskip

\noindent{\it Proof of Proposition \ref{geod}.} Let first $D$ be visible. The uniform equicontinuity
of the family $\F$ follows by Proposition \ref{equi}. Assume that $\bd D$ contains a non-trivial affine disc
$\DD.$ Then $to+(1-t)\DD\in D,$ where $o\in D$ and $0<t<1,$ which implies that
$\ds\liminf_{z\to p,w\to q}k_D(z,w)<\infty$ for any $p,q\in\DD.$ This contradicts \eqref{grom}.
\medskip

Let now the family $\F$ be uniformly equicontinuous and $\bd D$ be disc-free. Assume that $D$ is not visible.
Then we may find sequences $z_j\to p$ and $w_j\to q$ such that $(z_j|w_j)_o\to\infty,$
$p\neq q\in\bd D,$ $o\in D.$ Take $\vp_j\in \F$ with $\vp_j(\z_j)=z_j$ and $\vp_j(\e_j)=w_j.$
Passing to a subsequence, we may suppose $\vp_j\to\vp$ uniformly on $\D,$ where either
$\vp\in\F,$ or $\vp\in\O(\D,\bd D)$ The last case is impossible, since $p\neq q\in\vp(\ov\D)$
and $\bd D$ is disc-free. Then $\z_j\to\z,$ $\e_j\to\e,$ $\z\neq\e\in\bd\D,$ $o_j=\vp_j(0)\to o\in D,$
and $$(\z_j|\e_j)_0=(z_j|w_j)_{o_j}\to\infty$$
which contradicts the visibility of $\D.$\qed
\smallskip

\noindent{\it Remark.} The proof above shows that Proposition \ref{geod} remains true for any
bounded complete domain $D$ with Lipschitz boundary (which allows to move discs from $\bd D$ to $D$).
\medskip

\noindent{\it Proof of Proposition \ref{convex}.} First, we will show that
\eqref{quo} holds if $z_j,w_j\to p.$

We may take a complex geodesic $\vp_j$ through $z_j$ and $w_j$ such that $\d_D(a_j)=\max(\d_D\circ\vp_j),$
where $a_j=\vp_j(0).$

Then (the proof of) Proposition \ref{geod} implies that
$$\mbox{either }\d_D(a_j)\ge c_1>0,\mbox{ or }\vp_j\to p\mbox{ uniformly on }\D.$$

Since $\vp_j(D)\subset D\cap U$ for $j\gg 1$ in the second case, it follows that
$k_{D\cap U}(z_j,w_j)=k_D(z_j,w_j).$

Now we will consider the first case. It follows by \eqref{low} that
$$\frac{1}{2}\log\frac{c_1}{\d_D(\vp_j(\z))}\le k_D(a_j,\vp_j(\z))=k_\D(0,\z)<\frac{1}{2}\log\frac{2}{\d_\D(\z)}$$
and then, setting $c_2=2/c_1,$
\begin{equation}\label{1}
\d_\D(\z)\le c_2\d_D(\vp_j(\z)).
\end{equation}

On the other hand, for $\z_j=\vp_j^{-1}(z_j)$ and $\e_j=\vp_j^{-1}(w_j),$ the inequality
\begin{equation}\label{n1}
\log\left(1+\frac{|\z_j-\e_j|}{2\d_\D(\z_j)^{1/2}\d_\D(\e_j)^{1/2}}\right)\le k_\D(\z_j,\e_j).
\end{equation}
(see e.g. \cite[Lemma 4 (b)]{NTr}) and \eqref{dini} imply that
\begin{equation}\label{2}
\frac{|\z_j-\e_j|}{\d_\D(\z_j)^{1/2}\d_\D(\e_j)^{1/2}}\le\frac{c_3|z_j-w_j|}
{\d_D(z_j)^{1/2}\d_D(w_j)^{1/2}}
\end{equation}

Then, by \eqref{1},
\begin{equation}\label{3}
|\z_j-\e_j|\le c_2c_3|z_j-w_j|.
\end{equation}

Observe now that $\z_j\to\bd\D$ since $D$ is complete. After rotations we may assume that $\z_j\to 1$
and hence $\e_j\to 1$ by \eqref{3}. Proposition \ref{geod} provides $r>0$ such that
$\vp_j(\DD_r)\subset D\cap U,$ where $\DD_r=\D\cap\D(1,r)$ and $j\gg 1.$ Then
\begin{equation}\label{4}
\frac{k_{D\cap U}(z_j,w_j)}{k_D(z_j,w_j)}\le\frac{k_{\DD_r}(\z_j,\e_j)}{k_\D(\z_j,\e_j)},
\end{equation}
\begin{equation}\label{new}
\frac{k_{\DD_r}(\z_j,\e_j)}{k_\D(\z_j,\e_j)}
\le 1+c_4f_\D(\z_j,\e_j)
\end{equation}
for $j\gg 1.$ The proof of \eqref{new} is given at the end of this proof.

Note that $h$ is an increasing function. Since $h(c_3x)\le c_5h(x),$ then
\eqref{1}, \eqref{2}, \eqref{4}, and\eqref{new} imply that
$$\frac{k_{D\cap U}(z_j,w_j)}{k_D(z_j,w_j)}\le 1+c_2^2c_4c_5f_D(z_j,w_j).$$
So \eqref{quo} is proved when $z_j,w_j\to p.$

The proof of the first inequality in \eqref{met} (the second one is trivial) in this
case, which we skip, is similar and even simpler, taking into account the inequality
$$\k_{\DD_r}(z;\vec{e})\le(1+c_6\d^2_\D(z))\k_\D(z;\vec{e}),\quad z\in\DD_{r/2}$$
which follows by mapping $\DD_r$ conformally onto $\D.$

We may also deduce \eqref{met} from \eqref{quo}, shrinking $U$ such that $D\cap U$ to be convex;
then $\k_D$ and $\k_{D\cap U}$ are the ``derivatives'' of $k_D$ and $k_{D\cap U}.$

Thus we may find a neighborhood $V_p\Subset U$ of $p$ such that \eqref{quo} and \eqref{met}
hold if $z,w\in V_p.$ The same is true for any $p'\in\bd D\cap U.$
Then \eqref{met} follows by a compactness argument. On the other hand, \eqref{dini}
(applied to a Dini-smooth domain $\Omega$ such that $D\cap\ov{V}\subset\Omega\subset D\cap U$)
and \eqref{vis} imply that
$$\frac{k_{D\cap U}(z,w)}{k_D(z,w)}\le 1+c_{p',p''}f_D(z,w)$$
for $z\in D$ near $p'$ and $w\in D$ near $p'',$ where $p',p''\in\ov{D\cap V},$ $p'\neq p''.$

As we already know, the same is true when $p'=p''\in\partial D\cap\ov V.$ It is also valid if
$p'=p''\in D\cap\ov V$ because $k_{D\cap U}\le c'|z-w|$ and $k_D(z,w)\ge c''|z-w|$ for $z,w\in D$
near $p',$ where $c',c''>0.$

Assuming now that \eqref{quo} is not true, we may find points $p',p''\in\ov{D\cap V},$ and
sequences $(z_j)\to p',$ $(w_j)\to p'',$ $(c_j)\to\infty$ such that
$$\frac{k_{D\cap U}(z_j,w_j)}{k_D(z_j,w_j)}\ge 1+c_jf_D(z_j,w_j)$$
which is a contradiction.

It remains to observe that \eqref{lip} is a consequence of \eqref{quo} and \eqref{dini}.

\noindent{\it Subproof of \eqref{new}.} Mapping conformally $\DD_r$ onto $\D,$ one may see
as in \cite[Proposition 4.1 (4.5)]{NT2} (where the case of the upper half plane and its
intersection with $\D$ is discussed) that there exists $C_r>0$ such that
\begin{equation}\label{n2}
k_{\DD_r}(z,w)-k_\D(z,w)\le C_rg_\D(z,w),\ z,w\in\DD_{r/2}.
\end{equation}
Combining this inequality with \eqref{n1}, we get
\begin{equation}\label{n3}
\frac{k_{\DD_r}(z,w)}{k_\D(z,w)}\le 1+2C_r f_\D(z,w),\quad z,w\in\DD_{r/2}.
\end{equation}
\qed
\medskip

\noindent{\it Remarks.} (a) It follows as in (the proof) of \cite[Proposition 4.1]{NT2}
that an opposite inequality to \eqref{n2} holds if $z,w\in(1-r/2,1).$ Then \eqref{dini}
implies an opposite inequality to \eqref{n3} for such points, too.

(b) Call a point $p\in\bd D$ visible if \eqref{grom} is true for this $p$
with any point $q\in\bd D,$ $q\neq p.$
The above proof can be easily adapted to get the following version of Proposition~\ref{convex}
near a Dini-smooth visible point:
\smallskip

{\it Let $D$ be a bounded convex domain with Dini-smooth boundary near a visible point
$p\in\bd D.$ Then for any neighborhood $U$ of $p$ there exist $C>0$ and a neighborhood $V\Subset U$
of $p$ such that \eqref{quo}, \eqref{lip}, and \eqref{met} hold.}
\smallskip

For example, $p$ is a visible point if $\bd D\cap T^\C_p(\bd D)=\{p\}$ (see \cite[Lemma 4.5]{Zim}),
or \eqref{gold} is true near $p$ (see the proof of \cite[Theorem 1.4]{BZ1}).
\smallskip

(c) If $D'$ is a bounded domain with $\CC^{1+\eps}$ boundary, then
$$l_{D'}(z,w)\le\frac{1}{2}\log\frac{1}{\d_{D'}(z)}+\frac{1}{2}\log\frac{1}{\d_{D'}(w)}+C$$
by \cite[Theorem 1]{NPT}. The proof there can be modified to get the same if
$\partial D'$ is Dini-smooth (replacing the model domain $G^i$ defined on \cite[p.~427]{NPT}
by the model domain $G$ defined on \cite[p.~48]{NA}). Then the proofs of \eqref{quo} and \eqref{lip}
imply that they hold with $l_{D\cap U}$ instead of $k_{D\cap U}.$
\medskip

\noindent{\it Proof of Proposition \ref{local}.} By \cite[Proposition 1]{For} (when $\ov D$
is a Stein compact) and \cite[Lemma 5]{H1} (when $\bd D$ is $\CC^\infty$-smooth),
there exist a neighborhood $U_0$ of $p,$ a holomorphic map $\Phi:D\cup U_0\to\C^n,$ and a strongly
convex domain $D'\supset\Phi(D)$ such that $\Phi|_{U_0}$ is 1-1 and $\Phi(D\cap U_0)=D'\cap\Phi(U_0).$ Set
$u'=\Phi(u).$ Since $c_D(z,w)\le c_{D'}(z',w'),$
$k_{D\cap U}(z,w)=k_{D'\cap\Phi(U)}(z',w'),$ $\d_D(z)\sim\d_{D'}(z'),$ and $|z-w|\sim |z'-w'|$
if $z,w\in V\Subset U\subset U_0,$ we may assume that $D$ is strongly convex.
Then Proposition \ref{convex} is applicable.\qed
\medskip

\noindent{\it Proof of Theorem \ref{global}.} Assume that \eqref{two} is not true. Then we
may find sequences $(z_j)\to p\in\ov D,$ $(X_j)\subset C^n$ and $(c_j)\to\infty$ such that
$$\k_D(z_j;X_j)\ge(1+C\d^2_D(z))\g_D(z_j;X_j).$$
If $p\in D,$ then $\k_D(z;X)\le c'|X|$ and $\g_D(z;X)\ge c''|X|$ for $z,w\in D$
near $p',$ where $c',c''>0,$ which is a contradiction. If $p\in\partial D,$ we get
a contradiction to \eqref{met}.

Assume that \eqref{zero} is not true. Then we may find sequences $(z_j)\to p'\in\ov D,$
$(w_j)\to p''\in\ov D$ and $(c_j)\to\infty$ such that
$$\frac{l_D(z_j,w_j)}{c_D(z_j,w_j)}\ge 1+c_jf_D(z_j,w_j).$$
This is impossible if $p'\neq p'',$ since $l_D-c_D$ is bounded by \eqref{root}, and hence
$\frac{l_D(z,w)}{c_D(z,w)}$ is bounded for $z\in D$ near $p'$ and $w\in D$ near $p''.$
If $p'=p''\in D,$ the same is true, since $l_D(z,w)\le\tilde c|z-w|$ and $l_D(z,w)\ge\hat c|z-w|$
for $z,w\in D$ near $p',$ where $\tilde c,\hat c>0.$ When $p'=p''\in\partial D,$ we get a contradiction
to Proposition \ref{local} \eqref{quo}, applying to $U$ such that $\Phi(U)$ is convex
(then $l_D\le l_{D\cap U}=k_{D\cap U}$).

It remains to observe that \eqref{one} is a consequence of \eqref{zero} and \eqref{dini}.
\medskip

\noindent{\bf Acknowledgments.} The author would like to thank Pascal J. Thomas for useful discussions,
and the referees for their remarks which helped improve the exposition.

\end{document}